\documentclass{article}

\usepackage{amsfonts}
\usepackage{amssymb}
\usepackage{enumerate}
\usepackage{amsmath}
\usepackage[T1]{fontenc}
\usepackage{graphicx}

\newtheorem{ttt}{Theorem}[section]
\newtheorem{llll}[ttt]{Lemma}
\newtheorem{ccc}[ttt]{Claim}
\newtheorem{eee}[ttt]{Example}
\newtheorem{fff}[ttt]{Fact}
\newtheorem{rrr}[ttt]{Remark}
\newtheorem{sss}[ttt]{Statement}
\newtheorem{ddd}[ttt]{Definition}
\newtheorem{qqq}[ttt]{Question}
\newtheorem{cccc}[ttt]{Corollary}
\newtheorem{nnn}[ttt]{Notation}
\newtheorem{ccccc}[ttt]{Conjecture}

\newcommand{\beq}{\begin{equation} }

\newcommand{\bt}{\begin{ttt}}
\newcommand{\bl}{\begin{llll}}
\newcommand{\bc}{\begin{ccc}}
\newcommand{\bex}{\begin{eee}}
\newcommand{\bfa}{\begin{fff}}
\newcommand{\br}{\begin{rrr}\upshape}
\newcommand{\bst}{\begin{sss}}
\newcommand{\bd}{\begin{ddd}\upshape}
\newcommand{\bq}{\begin{qqq}}
\newcommand{\bnn}{\begin{nnn}}
\newcommand{\bcor}{\begin{cccc}}
\newcommand{\bcon}{\begin{ccccc}}

\newcommand{\eeq}{\end{equation}}
\newcommand{\et}{\end{ttt}}
\newcommand{\el}{\end{llll}}
\newcommand{\ec}{\end{ccc}}
\newcommand{\eex}{\end{eee}}
\newcommand{\efa}{\end{fff}}
\newcommand{\er}{\end{rrr}}
\newcommand{\est}{\end{sss}}
\newcommand{\ed}{\end{ddd}}
\newcommand{\eq}{\end{qqq}}
\newcommand{\ecor}{\end{cccc}}
\newcommand{\econ}{\end{ccccc}}
\newcommand{\enn}{\end{nnn}}

\newcommand{\bp}{\noindent\textbf{Proof. }}
\newcommand{\ep}{\hspace{\stretch{1}}$\square$\medskip}

\newcommand{\lab}[1]{\label{#1}}

\newcommand{\NN}{\mathbb{N}}
\newcommand{\ZZ}{\mathbb{Z}}
\newcommand{\QQ}{\mathbb{Q}}
\newcommand{\RR}{\mathbb{R}}

\newcommand{\al}{\alpha}

\newcommand{\e}{\varepsilon}

\newcommand{\si}{\sigma}

\newcommand{\iG}{\mathcal{G}}
\newcommand{\iH}{\mathcal{H}}

\newcommand{\sm}{\setminus}
\newcommand{\eset}{\emptyset}

\newcommand{\beeq}{\begin{equation}}
\newcommand{\eeeq}{\end{equation}}
\def\su{\subset}

\numberwithin{equation}{section}

\title{The structure of rigid functions}

\author{Rich\'ard Balka \\
\\
E\"otv\"os Lor\'and University\\
Department of Analysis\\
P\'az\-m\'any P. s. 1/c, H-1117, Budapest, Hungary\\
balkar@cs.elte.hu \\
\\
and\\
\\
M\'arton Elekes\thanks{Partially supported by
    Hungarian Scientific Foundation grants no.~49786 and 61600.}\\
\\
Alfr\'ed R\'enyi Institute of Mathematics\\
Hungarian Academy of Sciences\\
P.O. Box 127, H-1364 Budapest, Hungary\\
emarci@renyi.hu\\
www.renyi.hu/ $\tilde{}$ emarci}

\begin{document}

\maketitle

\begin{abstract}
A function $f:\RR \to \RR$ is called \emph{vertically rigid} if
$graph(cf)$ is isometric to $graph (f)$ for all $c \neq 0$. We
prove Jankovi\'c's conjecture by showing that a continuous
function is vertically rigid if and only if it is of the form
$a+bx$ or $a+be^{kx}$ ($a,b,k \in \RR$). We answer a question of
Cain, Clark and Rose by showing that there exists a Borel
measurable vertically rigid function which is not of the above
form. We discuss the Lebesgue and Baire measurable case,
consider functions bounded on some interval and
functions with at least one point of continuity. We also
introduce horizontally rigid functions, and show that a certain
structure theorem can be proved without assuming any regularity.
\end{abstract}

\insert\footins{\footnotesize{MSC codes: Primary 26A99 Secondary 39B22, 39B72,
    51M99}}
\insert\footins{\footnotesize{Key Words: rigid, functional equation,
    transformation, exponential, Lebesgue, Baire, Borel, ergodic}}

\section{Introduction}

An easy calculation shows that the exponential function $f(x) = e^x$ has the
somewhat `paradoxical' property that $cf(x)$ is a translate of $f(x)$ for
every $c>0$. It is also easy to see that every function of the form
$a+be^{kx}$  has this property.
This connection is also of interest from the point of view of
functional equations. In \cite{CCR} Cain, Clark and Rose introduced the notion
of vertical rigidity as follows.

\bd
A function $f:\RR \to \RR$ is called \emph{vertically rigid}, if $graph(cf)$ is
isometric to $graph (f)$ for all $c \in (0,\infty)$. (Clearly, $c \in \RR \sm
\{0\}$ would be the same.)
\ed

Obviously every function of the form $a+bx$ is also vertically
rigid. D.~Jankovi\'c conjectured (see \cite{CCR}) that the
converse is also true for continuous functions.

\bcon
\textbf(D.~Jankovi\'c)
A continuous function is vertically rigid if and
only if it is of the form $a+bx$ or $a+be^{kx}$ ($a,b,k \in \RR$).
\econ

The main result of the present paper is the proof of this conjecture.

We will need the following technical
generalisations.

\bd
If $C$ is a subset of $(0, \infty)$ and $\iG$ is a set of isometries of the
plane then we say that $f$ is vertically rigid \emph{for a set $C$ via
elements of $\iG$} if for every $c \in C$ there exists a $\varphi \in \iG$
such that $\varphi(graph(cf)) = graph (f)$.

(If we do not mention $C$ or $\iG$ then $C$ is $(0,\infty)$ and $\iG$
is the set of all isometries.)
\ed


The paper is organised as follows. In Section \ref{s:Jan} we prove
Jankovi\'c's conjecture, even if we only assume that $f$ is a
continuous vertically rigid function for an uncountable set $C$. We
show that it is sufficient to assume that $f$ has at least one point
of continuity, provided that it is vertically rigid for $C$
\emph{via translations}. We also show that it is sufficient to
assume that $f$ is bounded on some nondegenerate interval, provided
that it is vertically rigid \emph{via translations} and
$C=(0,\infty)$. In Section \ref{s:Borel} we show that Jankovi\'c's
conjecture fails for Borel measurable functions. Our example also
answers a question  from \cite{CCR} that asks whether every
vertically rigid function is of the form $a+bx$ $(a,b \in \RR)$ or
$a+be^{g}$ for some $a,b \in \RR$ and additive function $g$. In
Section \ref{s:meas} we prove that every Lebesgue (Baire) measurable
function that is vertically rigid \emph{via translations} is of the
form $a+be^{kx}$ \emph{almost everywhere (on a comeagre set)}. The
case of general isometries remains open. We also prove that in many
situations the exceptional set can be removed. In Section
\ref{s:sets} we define the notion of a rigid set, discuss how it is
connected to the notion of a rigid function, and prove an ergodic
theory type result. In Section \ref{s:horiz} we define horizontally
rigid functions, and give a simple characterisation of those
functions that are horizontally rigid \emph{via translations}.
Finally, in Section \ref{s:open} we collect the open questions.

\section{Proof of Jankovi\'c's conjecture}
\lab{s:Jan}

\bt
\lab{t:Jan}
(\textbf{Jankovi\'c's conjecture}) A continuous function is vertically rigid
if and only if it is of the form $a+bx$ or $a+be^{kx}$ $(a,b,k \in \RR)$.
\et

\br
In fact, our proof will show that it is sufficient if $f$ is a continuous
function that is vertically rigid for some uncountable set $C$.
\er

It is of course very easy to see that these functions are vertically rigid and
continuous. The proof of the difficult direction goes through three theorems,
which are interesting in their own right. First we reduce the general case to
translations, then the case of translations to horizontal translations, and
finally we describe the continuous functions that are vertically rigid via
horizontal translations.

\bt
\label{t:r->t} Let $f:\RR \rightarrow \RR$ be a continuous
function vertically rigid for an uncountable set $C \subset
(0,\infty)$. Then $f$ is of the form $a+bx$ for some $a,b \in \RR$
or $f$ is vertically rigid for an uncountable set $D\subset
(0,\infty)$ via translations.
\et

\bp
Let $\varphi_c$ be the isometry belonging to $c \in C$.
First we show that we may assume that these isometries are
orientation preserving. If uncountably many of the $\varphi_c$'s are orientation
preserving then we are done by shrinking $C$.
Otherwise let $C'\subset C$ be uncountable so
that $\varphi_{c'}$ is orientation reversing for every $c' \in C'$. Fix $c_0'
\in C'$, then one can easily check that
$c_0'f$ is vertically rigid via orientation preserving
isometries for $C'' = \left\{ \frac{c'}{c_0'} : c'\in C'\right \}$.
Suppose that we have already proved the theorem in case all isometries are
orientation preserving. Then either $c_0'f$ is of the form $a+bx$, and then so
is $f$, or $c_0'f$ is vertically rigid for an uncountable set $D$ via
translations, but then so is $f$ itself (for the same set $D$, but
possibly different translations).

For a function $f$ let $S_f$ be the set of directions between pairs of
points on the graph of $f$, that is,

\[
S_f = \left\{ \frac{p-q}{|p-q|} : p,q \in graph(f),\  p \neq q \right\}.
\]

Clearly $S_f$ is a symmetric (about the origin) subset of the unit circle $S^1
\su \RR^2$. As $f$ is a function, $(0, \pm 1) \notin S_f$.
Since $f$ is continuous, it is easy to see that $S_f$ actually consists of
two (possibly degenerate) nonempty intervals. (Indeed, if $p = (x,f(x))$ and
$q = (y,f(y))$ then $x<y$ and $x>y$ define two connected sets, open
half planes in $\RR^2$, whose continuous images form $S_f$.)

An orientation preserving isometry $\varphi$ of the plane is either a
translation or a rotation. Denote by $ang(\varphi)$ the angle of $\varphi$ in
case it is a rotation, and set $ang(\varphi) = 0$ if $\varphi$ is a
translation.

Now we define two self-maps of $S^1$. Denote by $\varrho_\al$ the rotation about
the origin by angle $\al$. For $c>0$ let $\psi_c$ be the map obtained by
`multiplying by $c$', that is, let
\[
\psi_c((x,y)) = \frac{(x,cy)}{|(x,cy)|}\ \ \ ((x,y) \in S^1).
\]

It is easy to see that the rigidity of $f$ implies that for every $c \in C$
\beq
\lab{e:inv}
S_f = \varrho_{ang(\varphi_c)} (\psi_c(S_f)).
\eeq

If $S_f$ consists of two points, then $f$ is clearly of the form $a+bx$ and we
are done.

Let now $S_f = I \cup -I$, where $I$ is a subinterval of $S^1$ in
the right half plane. We claim that the endpoints of $I$ are among
$(0,\pm 1)$ and $(1,0)$. Suppose this fails, and consider
the function $l(c) =  arclength \left( \psi_c(I) \right)$ $(c \in
(0,\infty))$. It is easy to see that $l$ is real analytic,
and we show that it is not constant. Let us first assume that $(0,1)$
and $(0,-1)$ are not endpoints of $I$, then $\lim_{c \to 0} l(c) = 0$,
so $l$ cannot be constant (as $l>0$). Let us now suppose that either
$(0,1)$ or $(0,-1)$ is an endpoint of $I$, then $0 < arclength(I) <
\frac\pi2$
or $\frac\pi2 < arclength(I) < \pi$. In both cases $\lim_{c \to 0}
l(c) = \frac\pi2$ but $l(c) \neq \frac\pi2$, so $l$ is not
constant. As $l$ is analytic, it  attains each of its values at most
countably many
times, so there exists a  $c \in C$ so that $arclength
\left( \psi_c(I) \right) \neq arclength(I)$, which contradicts
(\ref{e:inv}).

(Actually, it can be shown by a somewhat lengthy calculation using the
derivatives that $l$ attains each value at most twice.)

But this easily yields $ang(\varphi_c) = 0$ or $\pi$ for every $c
\in C$. (Note that $(0,\pm 1) \notin  \nolinebreak S_f$ and that
$S_f$ is symmetric.) Just as above, we may assume that
$ang(\varphi_c) = 0$ for all $c \in C$. (Indeed, choose $C',c_0'$
analogously.) But then $f$ is vertically rigid for an uncountable
set via translations, so the proof is complete.
\ep

\bt
\label{t:t->ht}
Let $f:\RR \rightarrow \RR$ be an arbitrary
function that is vertically rigid for a set $C \su (0,\infty)$ via
translations. Then there exists $a\in \RR$ such that $f-a$ is
vertically rigid for the same set via horizontal translations.
\et

\bp
We can clearly assume that $1 \notin C$.
By assumption, for every $c \in C$ there exists $u_{c},v_{c} \in \RR$ such that

\beeq
\label{14}
cf(x)=f(x+u_{c})+v_{c} \ (\forall x \in \RR)
\eeq

Applying this first with $c=c_2$ then with $c=c_1$ we obtain

\[c_1c_2 f(x) =
c_1(f(x+u_{c_{2}})+v_{c_{2}}) =
\]

\beeq
\lab{17}
= c_{1}f(x+u_{c_{2}})+c_{1}v_{c_{2}} =
f(x+u_{c_{1}}+u_{c_{2}})+v_{c_{1}}+c_{1}v_{c_{2}}
\eeq

Interchanging $c_1$ and $c_2$ we get

\beeq
\label{18}
c_{2}c_{1} f(x) = f(x+u_{c_{2}}+u_{c_{1}})+v_{c_{2}}+c_{2}v_{c_{1}}.
\eeq

Comparing (\ref{17}) and (\ref{18}) yields
$
v_{c_{1}}+c_{1}v_{c_{2}}=v_{c_{2}}+c_{2}v_{c_{1}},
$
so

\[
\frac{v_{c_{1}}}{c_{1}-1}=\frac{v_{c_{2}}}{c_2-1}\ \ \
\textrm{ for all } c_{1},c_{2} \in C,
\]
consequently
$
a:=\frac{v_{c}}{c-1}
$
is the same value for all $c \in C$.
Substituting this back to (\ref{14}) gives
$cf(x) = f(x+u_{c})+a(c-1)$,
so $c(f(x)-a) = f(x+u_{c})-a$ for all $c\in C$,
%
hence $f-a$ is vertically rigid for $C$ via horizontal translations.
\ep

\bt
\label{t:htcont}
Let $f:\RR \rightarrow \RR$ be a continuous
vertically rigid function for an uncountable set $C\subset
(0,\infty)$ via horizontal translations. Then $f$ is of the form
$be^{kx}$ $(b\in \RR, k \in \RR\setminus\{0\})$.
\et

Before proving this theorem we need a definition and a lemma.

\bd
For a function $f:\RR\to\RR$ let $T_{f,C} \su \RR$ be the
additive group generated by the set $T' = \{ t\in\RR : \exists c \in C \
\forall x \in \RR \  f(x+t)=cf(x) \}$. (We will usually simply write $T$
for $T_{f,C}$.)
\ed

\bl
\lab{l:add}
Let $f:\RR \rightarrow \RR$ be a vertically rigid function for an uncountable
set $C\subset (0,\infty)$ via horizontal translations such that $f(0)=1$.
Then $T$ is dense and
\[
f(x+t) = f(x) f(t) \ \ \forall x \in \RR \ \forall t \in T.
\]
Moreover,
$f(t)>0$ for every $t \in T$.
%
\el

\bp
By assumption, for every $c \in C$ there exists $t_c \in \RR$
such that $cf(x) = f(x+t_c)$ for every $x \in \RR$. Then $t_c\in
T$ for every $c \in C$.  Since $f$ is not
identically zero, $t_c \neq t_{c'}$ whenever $c,c' \in C$
are distinct. Hence $\{ t_c : c \in C \}$ is uncountable, so $T$
is uncountable. As every subgroup of $\RR$ is either discrete
countable or dense, $T$ is dense.

Every $t\in T$ can be written as $t=\sum _{i=1} ^{m} n_{i}t_{i}$
($t_i\in T', n_i\in \ZZ, i=1,\dots ,m$) where
$f(x+t_{i})=c_{i}f(x) \  (x \in \RR, \ i=1,\dots ,m$).

From these we easily get \beeq \lab{2}  f(x+t)=c_{t}f(x),
\ \ \textrm{ where } c_t = \prod _{i=1} ^{m} c_{i}^{n_{i}}, \ x \in
\RR, \ t \in T.
\eeq
 Note that $c_{t}>0$
(and also that it is not necessarily a member of $C$). It suffices
to show that $c_t = f(t)$ for every $t \in T$, but this follows if we
substitute $x=0$ into (\ref{2}).
\ep

\bp (Thm. \ref{t:htcont}) If $f$ is identically zero then we are
done, so let us assume that this is not the case.  The class of
continuous vertically rigid functions for some uncountable set via
horizontal translations, as well as the class of functions of the
form $be^{kx}$ ($b\in \RR, k \in \RR\setminus\{0\}$) are both closed
under horizontal translations and under multiplication by nonzero
constants, so we may assume that $f(0)=1$. Then the previous lemma
yields that $f(t_1+t_2) = f(t_1)f(t_2)$ $(t_1,t_2 \in T)$, and also
that $f|_T>0$. Then $g(t) = \log f(t) $ is defined for every $t \in
T$, and $g$ is clearly additive on $T$. But it is well-known (and an
easy calculation) that an additive function on a dense subgroup is
either of the form $kx$, or unbounded both from above and below on
every nondegenerate interval. The second alternative cannot hold,
since $f$ is continuous, so $f|_T$ is of the form $e^{kx}$, so by
continuity $f$ is of this form everywhere. Since $C$ contains
elements different from $1$, we obtain that $f(x)=1$ $(x \in \RR)$
is not vertically rigid for $C$ via horizontal translations, hence
$k\neq 0$. \ep

Putting together the three above theorems completes the proof of Jankovi\'c's
conjecture.

We remark here that we have actually also proved the following,
which applies e.g.~to Baire class 1 functions.

\bt \label{c2} Let $f:\RR \rightarrow \RR$ be a vertically rigid
function for an uncountable set $C\subset (0,\infty)$ via
translations. If $f$ has a point of continuity then it is of the
form $a+be^{kx}$ ($a,b,k \in \RR$). If $f$ is vertically rigid via
translations (i.e.~$C=(0,\infty)$) and bounded on a nondegenerate
interval then it is of the form $a+be^{kx}$ ($a,b,k \in \RR$), too.
\et

\bp
Following the proof of the last theorem we may assume in both cases
that $f(0)=1$, the translations are horizontal, and $f|_T$ is of the
form $e^{kx}$ $(k \in \RR)$.

In the first case, let $x_{0}$ be a point of continuity of $f$, then
clearly $f(x_0)=e^{kx_0}$, since $T$ is dense. Let now $x\in \RR$ be
arbitrary, and $t_{n}\in T$ $(n\in \NN)$ be such that $\lim _{n \to
  \infty} t_{n} =x_{0}-x$. Using Lemma \ref{l:add} we obtain

\[
e^{kx_0} = f(x_0) = \lim _{n \to \infty} f(x+t_{n}) = \lim _{n \to
\infty} f(x)f(t_{n}) = f(x) \lim _{n \to \infty} f(t_{n}) =
\]

\[
= f(x) \lim _{n \to \infty} e^{kt_n} = f(x) e^{k(x_0-x)} = f(x) e^{kx_0}
/ e^{kx},
\]
from which $f(x) = e^{kx}$ follows.

In the second case, for every $c>0$ there is a $t_{c}\in
T=T_{f,(0,\infty)}$ such that $cf(x)=f(x+t_{c})=f(x)f(t_{c})$. By
substituting $x=0$ into the equation we get $c=f(t_{c})=e^{kt_{c}}$
for every $c>0$. (In particular, $k\neq 0$.) So $t_{c}=\frac{\log
c}{k}$. If $c$ ranges over $(0,\infty)$ then $t_{c}$ ranges over
$\RR$, so we get $T=\RR$. Hence $f|_T=f$ is of the form $e^{kx}$,
and we are done. \ep

\bex There exists a function $f:\RR \rightarrow \RR$ that is
vertically rigid for an uncountable set $C \subset \RR$ via
horizontal translations, bounded on every bounded interval, and is
\emph{not} of the form $a+be^{kx}$ ($a,b,k \in \RR$). \eex

\bp
Let $P\su\RR$ be an uncountable linearly independent set
over $\QQ$, see e.g. \cite[19.2]{Ke} or \cite{vN}.
Define $\widehat{P}$ to be the generated additive subgroup.
Let

\[
f(x) =
\begin{cases}

e^x & \textrm{ if } x\in \widehat{P}
\\ 0 & \textrm{ if } x\in \RR \setminus \widehat{P},
\end{cases}
\]
then $f$ is clearly bounded on every bounded interval.

It is easy to see that $\frac{p}2 \in \RR \setminus \widehat{P}$ for
every $p \in P$, so $\widehat{P} \neq \RR$, hence $f$ is not
continuous, so it is not of the form $a+be^{kx}$ ($a,b,k \in \RR$).

For every $p\in P$ and $x \in \RR$ we have
$x\in \widehat{P} \iff x+p\in \widehat{P}$, which easily implies
$f(x+p)=e^{p}f(x)$. Hence $f$ is vertically rigid for the
uncountable set $C=\{e^p: \, p\in P\}$.
\ep

Jankovi\'c's conjecture has the following curious corollary.

\bcor
There are continuous functions $f$ and $g$ with isometric graphs so that $f$
is vertically rigid but $g$ is not.
\ecor

\bp
If we rotate the graph of $f(x) = e^x$ clockwise by $\frac\pi4$ then we obtain
the graph of a continuous function. By Theorem \ref{t:Jan} it is not
vertically rigid.
\ep

\section{A Borel measurable counterexample}
\lab{s:Borel}

In this section we show that Jankovi\'c's conjecture fails for Borel
measurable functions. Our example also answers Question 1 in \cite{CCR} of
Cain, Clark and Rose, which asks whether every vertically rigid function is of
the form $a+bx$ $(a,b \in \RR)$ or $a+be^{g}$ for some $a,b \in \RR$ and
additive function $g$. By Thm. 2 of \cite{CCR} $a+be^{g}$ is vertically rigid
if and only if $b=0$ or $g$ is surjective.

\bt
\lab{t:count}
There exists a Borel measurable vertically rigid function $f: \RR \to
[0,\infty)$ (via horizontal translations) that is not of the form $a+bx$ $(a,b
\in \RR)$ or $a+be^{g}$ for some $a,b \in \RR$ and additive function $g$.
\et

For definitions and basic results on Baire measurable sets (= sets
with the property of Baire), meagre
(= first category) and comeagre (= residual) sets consult
e.g.~\cite{Ox} or \cite{Ke}. For Polish spaces and Borel
isomorphisms see e.g.~\cite{Ke}.

\bp
Let $P$ be a Cantor set (nonempty nowhere dense compact set
with no isolated points) that is linearly independent over $\QQ$, see
e.g.~\cite[19.2]{Ke}. (One can also derive the existence of such a set
from \cite{vN} using the well-known fact that every
uncountable Borel or analytic set contains a Cantor set.)
It is easy to see that for all $n_1, \dots, n_k \in \ZZ$
the set $ P_{n_1, \dots, n_k} = \{ n_1p_1 + \dots + n_kp_k : p_1,
\dots ,p_k \in P\}$ is compact, hence the group $\widehat{P}$
generated by $P$ (that is, the union of the $P_{n_1, \dots,
n_k}$'s) is a Borel, actually $F_\si$ set. As $P$ is linearly
independent, each element of $\widehat{P}$ can be uniquely written
in the form $n_1p_1 + \dots + n_kp_k$.

Since $P$ and $(0,\infty)$ are uncountable Polish spaces, we can
choose a Borel
isomorphism $g:P\rightarrow (0,\infty)$. Let $f:\RR \rightarrow \RR$
be the following function:

\[
f(x) =
\begin{cases}
0 & \textrm{ if } x\in \RR \setminus \widehat{P} \\
\prod _{i=1}^k g(p_i)^{n_i} & \textrm{ if } x = \sum_{i=1}
^k n_ip_i \in \widehat{P}, \ n_i \in \ZZ, \ p_i \in P,
\ i=1,\dots,k.
\end{cases}
\]

This function is Borel, as it is Borel on the countably many Borel sets
$P_{n_1, \dots, n_k}$, and zero on the rest. However, $f$ is not continuous,
as it is unbounded on the compact set $P$. Therefore $f$ is not of the form
$a+bx$. Suppose now that $f$ is of the form $a+be^{g}$ for some $a,b \in \RR$
and additive function $g$. Clearly $b \neq 0$, since $f$ is not constant,
therefore $\frac{f-a}{b} = e^g$ is Borel measurable, and then so is $g$ by
taking logarithm. But it is well-known that every Borel (or even Lebesgue)
measurable additive function is of form $kx$ $(k \in \RR)$, hence $f$ is
continuous, a contradiction.

What remains to show is that $f$ is vertically rigid via
horizontal translations. For every $c>0$ there exists a $p \in P$
such that $g(p)=c$. Now we check that $cf(x) = f(x+p)$ for all $x
\in \RR$. Clearly $x \in \widehat{P}$ if and only if $x+p \in
\widehat{P}$. Therefore $cf(x) = f(x+p) = 0$ if $x \notin
\widehat{P}$. Let now $x = n_1p_1 + \dots + n_kp_k \in
\widehat{P}$, and assume without loss of generality that $p=p_1$
($n_1 = 0$ is also allowed). Then $cf(x)  = g(p)f(x) =  g(p)
g(p)^{n_{1}} g(p_2)^{n_2} \cdots g(p_k)^{n_k} =  g(p)^{n_{1}+1}
g(p_2)^{n_2} \cdots g(p_k)^{n_k} =  f( (n_{1}+1)p + n_2p_2+\dots +
n_kp_k ) = f(x+p)$, which finishes the proof.
\ep

\section{Lebesgue and Baire measurable functions}
\lab{s:meas}

It is easy to see that the example in the previous section is zero almost
everywhere (on a comeagre set). Indeed, it can be shown that
every $P_{ n_1, \dots, n_k }$ has uncountably many pairwise disjoint
translates.

Therefore it is still possible that the complete analogue of
Jankovi\'c's conjecture holds: every vertically rigid Lebesgue
(Baire) measurable function is of the form $a+bx$ or $a+be^{kx}$
\emph{almost everywhere (on a comeagre set)}. In this section we
prove this in case of translations. The general case remains open,
see Section \nolinebreak \ref{s:open}. We also prove that in many
situations the exceptional set can be removed.

\bt
\label{t:meas}
Let $f:\RR \rightarrow \RR$ be a vertically rigid function
for an uncountable set $C\subset (0,\infty)$ \emph{via
translations}. If $f$ is Lebesgue (Baire) measurable then it is of the form
$a+be^{kx}$ $(a,b,k \in \RR)$ almost everywhere (on a comeagre set).
\et

\bp By Theorem \ref{t:t->ht} we can assume that $f$ is vertically
rigid for $C$ via horizontal translations. As in the proof of
Theorem \ref{t:htcont} we can also assume that $f(0) = 1$. Then
Lemma \ref{l:add} implies that \beq \lab{e:1} f(x+t) = f(x) f(t) \ \
\forall x \in \RR \ \forall t \in T \eeq and $f(t)>0$ for every $t
\in T$.

First we show that the sign of $f$ is constant almost everywhere
(on a comeagre set). It is easy to see from (\ref{e:1}) that the
sets $\{ f>0 \}$, $\{ f=0 \}$, and $\{ f<0 \}$ are all Lebesgue
(Baire) measurable sets periodic modulo every $t \in T$. It is a
well-known and easy consequence of the Lebesgue density theorem
(the fact that every set with the Baire property is open modulo
meagre) that if a measurable set $H$ has a dense set of periods
then either $H$ or $\RR \sm H$ is of measure zero (meagre). But
the above three sets cover $\RR$, hence at least one of them is of
positive measure (nonmeagre), and then that one is of full measure
(comeagre). If $f=0$ almost everywhere (on a comeagre set) then we
are done, otherwise we may assume that $f>0$ almost everywhere (on
a comeagre set). (Indeed, $-f$ is also rigid via horizontal
translations, and then we can apply a horizontal translation and a
positive multiplication to achieve $f(0) = 1$.)

Set $D = \{ f>0 \}$ and define the measurable function $g = \log f$
on $D$. Recall that $D+t = D$ ($\forall t \in T$) and note that $T
\su D$. Clearly

\[
g(x+t) = g(x) + g(t) \ \ \forall x \in D \  \forall t \in T,
\]
so $g|_T$ is additive. Now we show that $g|_T$ is of the form $kx$.
Let us suppose that this is not the case. As we have mentioned
above, if an additive function is not of the form $kx$ then it is
unbounded on every interval from above (and also below). For every
Lebesgue (Baire) measurable function there is a measurable set of
positive measure (nonmeagre) on which the function is bounded. So
let $M \su D$ be a measurable set of positive measure (nonmeagre)
such that $\left| g|_M \right| \le K$ for some $K \in \RR$. By the
Lebesgue density theorem (the fact that every Baire measurable set
is open modulo meagre) there exists $\e > 0$ so that $( M+s ) \cap M
\neq \eset$ for every $s \in (-\e,\e)$. Choose $t_0\in T$ in
$(-\e,\e)$ so that $g(t_0) > 2K$. Fix an arbitrary $m_0 \in M\cap
(M-t_{0})$, then $g(m_0+t_0) = g(m_0) + g(t_0) > g(m_0) + 2K$, which
is absurd, since $m_0+t_0, m_0 \in M$ and $\left| g|_M \right| \le
K$.

Now define $h(x) = g(x) - kx$ ($x \in D$). This is a measurable
function that is periodic modulo every $t \in T$. Indeed,
\[
h(x+t) = g(x+t) - k(x+t) = g(x) -kx + g(t) - kt = h(x) + 0 = h(x).
\]

It is a well-known consequence of the Lebesgue density
theorem (the fact that every Baire measurable set is open modulo meagre)
that if the periods of a measurable function form a dense set then
the function is constant almost everywhere (on a comeagre set). Hence $g(x) =
kx + c$ almost everywhere (on a comeagre set), so $f(x) = e^c e^{kx}$ almost
everywhere (on a comeagre set), so we are done.
\ep

Our next theorem shows that the measure zero (meagre) set can be
removed, unless $f$ is constant almost everywhere (on a comeagre
set). Theorem \nolinebreak \ref{t:count} provides an almost
everywhere (on a comeagre set) constant but nonconstant function
that is vertically rigid via horizontal translations.

\bt \lab{xy} Let $f:\RR \to\RR$ be a vertically rigid function that
is of the form $a+bx$ ($b \neq 0$) or $a+be^{kx}$ ($bk \neq 0$)
almost everywhere (on a comeagre set). Then $f$ is of this form
everywhere. \et

Let us denote the 1-dimensional Hausdorff measure by $\iH^1$. For the
definition and properties see \cite{Fa} or \cite{Ma}. First we
prove the following lemma.

\bl \label{lem} Let $f, g : \RR \to \RR$ be arbitrary functions, and
let $\varphi$ be an isometry such that $\varphi( graph( f )) =
graph( g )$.  Let $f', g' : \RR\rightarrow \RR$ be continuous
functions such that $f' = f$ almost everywhere (on a comeagre set)
and $g' = g$ almost everywhere (on a comeagre set). Let us also
assume that $graph(f')$, $\varphi( graph(f'))$, $graph(g')$ and
$\varphi^{-1}( graph(g'))$ are coverable by the graphs of countably
many Lipschitz (continuity suffices for the category case)
functions. Then $\varphi( graph(f')) = graph(g')$. \el

\bp
By symmetry of $f'$ and $g'$ (with $\varphi^{-1}$), it suffices to
show that $graph(g') \su \varphi( graph( f' ) )$. Since the latter set
is closed, it also suffices to show that $\varphi(
graph( f' ))$ covers a dense subset of $graph(g')$. We will actually
show that $\varphi( graph( f' ))$ covers $\iH^1$ a.e.~(relatively
comeagre many) points of $graph(g')$, which will finish the proof.

If an element of $graph(g')$ fails to be covered by $\varphi(graph(
f' ))$ then it is either in $graph(g') \sm graph(g)$ or in $\varphi(
graph(f) \sm graph(f'))\cap graph(g')$. The first set is clearly of
$\iH^1$ measure zero (relatively meagre in $graph(g')$), so it
suffices to show that this is also true for the second.
Equivalently, we need that $graph(f) \sm graph(f')$ only covers a
$\iH^1$ measure zero (relatively meagre) subset of $\varphi^{-1}(
graph(g') )$. Suppose that $\varphi^{-1}( graph(g') ) \su
\bigcup_{n=1}^\infty graph(h_n)$, where the $h_n$'s are Lipschitz
(continuous) functions. As $graph(h_n) \cap (graph(f) \sm
graph(f'))$ is clearly of $\iH^1$ measure zero for every $n$, we are
done in the measure case.

Let us now write $\{x \in \RR: f'(x) \neq f(x) \} =
\bigcup_{m=1}^\infty N_m$, where each $N_m$ is nowhere dense. It is
enough to show that each $graph(f|_{N_m})$ only covers a relatively
nowhere dense subset of $\varphi^{-1}( graph(g'))$. Fix an $m$, and
suppose that $graph(f|_{N_m})$ is  dense in an open subarc $U \su
\varphi^{-1}( graph(g'))$. By the Baire Category Theorem there
exists a relatively open subarc $V \su U$ that is covered by one of
the $graph(h_n)$'s. But this is impossible, as the arc $V$ is in
$graph(h_n)$, and the set $N_m \su \RR$ is nowhere dense, so even
$N_m \times \RR$ covers at most a relatively nowhere dense subset of
$V$, hence $graph(f|_{N_m})$ cannot be dense in $V$. \ep

\bp (Thm. \ref{xy}) Using the notation of the above lemma, let first
$f$ be a vertically rigid function such that $f=f'$ almost
everywhere (on a comeagre set), where $f'$ is of the form 
$a+be^{kx}$ ($bk \neq 0$). 
The above lemma implies that $f'$ is also vertically rigid with
the same isometries $\varphi_c$. By considering the unique asymptote
and the limit at $\pm \infty$ of $f'$ we obtain that every
$\varphi_c$ is a translation. By Theorem \ref{t:t->ht} we may assume
that every $\varphi_c$ is actually horizontal, hence $f'$ is of the form 
$be^{kx}$. Hence $cf'(x) = f'\left(x+\frac{\log(c)}{k}\right)$ for
every $x \in \RR$, $c>0$ and the same holds for $f$. Assume now that
there is an $x_0$ so that $f(x_0) \neq f'(x_0)$, then $cf(x_0) \neq
cf'(x_0)$ for every $c>0$, therefore
$f\left(x_0+\frac{\log(c)}{k}\right) \neq
f'\left(x_0+\frac{\log(c)}{k}\right)$ for every $c>0$, which is a
contradiction as $f = f'$ almost everywhere (on a comeagre set).

Assume now that $f'$ is of the form $a+bx$ ($b \neq 0$). 
First we show that $f'$
is vertically rigid by the same isometries as $f$. For every $c>0$
set $g = cf$, $g' = cf'$, and let $\varphi_c$ be the isometry
mapping $graph(f)$ onto $graph(g)$. As $graph(f) \cap graph(f')$
contains at least two points and $\varphi_c(graph(f) \cap
graph(f'))$ is the graph of a function we obtain that the line
$\varphi_c(graph(f'))$ is not vertical, and similarly for
$\varphi_c^{-1}( graph(g'))$. Therefore they are coverable by the
graphs of countably many, actually a single, Lipschitz (continuous)
function, hence the previous lemma applies. Hence $f'$ is vertically
rigid by the same isometries as $f$.

Similarly to Theorem \ref {t:r->t} we can assume that $f$ is
vertically rigid via orientation preserving isometries for a set $C$
of positive outer measure (nonmeagre). So $\varphi_c$ is a rotation
or translation for every $c\in C$, and by splitting $C$ into two
parts and keeping one with positive outer measure (nonmeagre), we
can assume that $A=\{ang(\varphi _{c}): c\in C\}$ is a subset of the
left or the right half of the unit circle. We could calculate
$ang(\varphi _{c})$ explicitly, but we only need that it is a 
nonconstant real analytic function. From this it is easy to see that
the set $A$ is of positive outer measure (nonmeagre). Assume now
that there is an $x_0$ so that $f(x_0) \neq f'(x_0)$. We prove that
this contradicts the fact that $\varphi_c( graph(f) )$ is the graph
of a function for every $c\in C$. For this it suffices to show that
$S_f$ (see Theorem \ref {t:r->t}) is of full measure (comeagre). But
this clearly follows simply by looking at the pairs $(p_0,q)$ and
$(q,p_0)$ where $p_0 = (x_{0}, f(x_0))$ and $q$ ranges over
$graph(f) \cap graph(f')$. \ep

\section{Rigid sets}
\lab{s:sets}

The starting point is the proof of Theorem \ref{t:r->t}. So far we
are only able to prove this result for continuous functions, and
consequently we can only handle translations in the
Borel/Lebesgue/Baire measurable case. But generalisations of the
ideas concerning the sets $S_f$ could tackle this difficulty. For a
Borel function $f$ the set $S_f$ is analytic (see e.g.~\cite{Ke}),
and every analytic set has the Baire property, so the result of this
section can be considered as the first step towards handling Borel
functions with general isometries.
See Equation (\ref{e:inv})
for the following notations.

\bd
We call a symmetric (about the origin) set $H \subset S^1$
\emph{rigid} for a set $C \su (0,\infty)$ if for every $c\in C$
there is an $\al$ such that
\beq
\lab{e:3}
H = \varrho_\al
(\psi_c(H)).
\eeq
\ed

\bl
Let $U$ be a regular open set (i.e. $int(cl(U))=U$) that is
rigid for an uncountable set $C$. Then $U=\emptyset$, or $U=S^1$, or
every connected component of $U$ is an interval whose endpoints
are among $(0,\pm 1)$ and $(\pm 1,0)$.
%
\el

\bp
Let $A$ be the set of arclengths of the connected components of $U$,
then $A$ is countable. Let $I$ be a connected
component of $U$ showing that $U$ is not of the desired form, then
$0 < arclength(I) < \pi$ since $U$ is symmetric and regular. As in the
proof of Theorem \ref{t:r->t} let us prove that the real analytic
function $l(c) = arclength( \psi_c(I) )$ $(c \in (0,\infty))$ is not
constant. If $I$ is in the left or right half of $S^1$ then we already
showed this there, so we may assume that $(0,1)$ or $(0,-1)$ is in
$I$. Since $\lim_{c \to \infty} \psi_c(x) \in \{ (0,\pm 1) , (\pm
1,0) \}$ for every $x \in S^1$, we obtain that $\lim_{c \to \infty}
l(c) \in \ZZ \frac\pi2$. Hence we are done using $0 <
arclength(I) < \pi$ unless $arclength(I) = \frac\pi2$. But if
$arclength(I) = \frac\pi2$ then $\lim_{c \to \infty} l(c) = 0$
since $(0,1)$ or $(0,-1)$ is in $I$, and therefore $l$ cannot be constant.

Hence $l$ attains each of its values at most countably many times,
so there is a $c \in C$ such that $arclength( \psi_c(I) ) \notin A$,
contradicting (\ref{e:3}).
\ep

One can also show, using an argument similar to the above one (by
considering the possible
distances of pairs in $H$), that the rigid sets (for $C = (0,\infty)$)
of cardinality
smaller than the continuum are the following: the empty set, the
symmetric sets of two elements and the set $\{ (0,\pm 1) , (\pm
1,0) \}$.

The next statement is somewhat of ergodic theoretic flavour.

\bt
Let $H$ be a Baire measurable set that is
rigid for an uncountable set $C$. Then in each of the four
quarters of $S^1$ determined by $(0,\pm 1)$ and $(\pm 1,0)$
either $H$ or $S^1 \sm H$ is meagre.
\et

\bp
$H$ can be written as $H = U \Delta F$ in a unique way, where $U$ is regular
open, $F$ is meagre and $\Delta$ stands for symmetric difference,
see \cite[4.6]{Ox}. Then it is easy to see by the
uniqueness of $U$ that $U$ is rigid for $C$, so we are done by the
previous lemma.
\ep

\section{Horizontally rigid functions}
\lab{s:horiz}

In this section we characterise the functions that are horizontally rigid via
translations. This answers Question 3 of \cite{CCR} in the case of
translations.

\bd
A function $f:\RR \rightarrow \RR$ is \emph{horizontally rigid}, if
$graph(f(cx))$ is isometric to $graph(f(x))$ for all $c \in (0,\infty)$.
\ed

\bt
\lab{t:horiz}
A function $f:\RR \rightarrow \RR$ is
horizontally rigid via translations if and only if there exists $r
\in \RR$ such that $f$ is constant on $(-\infty, r) $ and
$(r,\infty)$.
\et

\bp
These functions are trivially horizontally rigid via translations. As the
proof of the other direction resembles that of
Theorem \ref{t:t->ht}, we only sketch it.

For every $c>0$ there exist $u_{c},v_{c} \in \RR$ such that
$f(cx)=f(x+u_{c})+v_{c} \linebreak \ \ (x \in  \RR)$.
We may assume $u_{1}=v_{1}=0$. If $c \in (0,\infty)
\setminus \{1\}$
then there is an $x_{c} \in  \RR$ such that $cx_{c}=x_{c}+u_{c}$, and
substituting this back to the above equation we get $v_{c}=0$. Hence
$f(cx)=f(x+u_{c}) \ \ (x \in \RR)$
for every $c \in (0,\infty)$.

First we show that if $f$ has a period $p>0$ then $f$ is constant.
Using the last equation twice we obtain
\[
f(cx)=f(x+u_{c})=f(x+u_{c}+p)
=f\left((x+p)+u_{c}\right) = f\left(c(x+p)\right) = f(cx+cp).
\]

If $x$ ranges over $\RR$ then so does $cx$, hence $cp$ is also a period. If
$c$ ranges over $(0,\infty)$, then so does $cp$, hence every positive
number is a period, so $f$ is constant.

Using $f(cx)=f(x+u_{c})$ again twice we obtain
\[
f(c_{1}(c_{2}x)) = f(c_{2}x+u_{c_{1}}) =
f\left(c_{2}\left(x+\frac{u_{c_{1}}}{c_{2}}\right)\right) =
f\left(x+\frac{u_{c_{1}}}{c_{2}}+u_{c_{2}}\right).
\]

Interchanging $c_1$ and $c_2$ and comparing the two equations we get
\[
f\left(x+\frac{u_{c_{1}}}{c_{2}}+u_{c_{2}}\right) =
f\left(x+\frac{u_{c_{2}}}{c_{1}}+u_{c_{1}}\right),
\]
so
$\pm
\left[
  \left(\frac{u_{c_{1}}}{c_{2}}+u_{c_{2}}\right) -
  \left(\frac{u_{c_{2}}}{c_{1}}+u_{c_{1}}\right)
\right]$
is a period, and hence it is zero. Therefore

\[
\frac{u_{c_{1}}}{1 - \frac{1}{c_{1}}} = \frac{u_{c_{2}}}{1- \frac{1}{c_{2}}}
\textrm{ for every } c_{1},c_{2} \in (0,\infty) \setminus \{1\}.
\]

Set $r = \frac{u_{c}}{1 - \frac{1}{c}}$, then
$u_{c} = r\left(1-\frac{1}{c}\right)$
for every $c \in (0,\infty)$. Substituting this back to $f(cx)=f(x+u_{c})$
gives
$f(cx) = f \left( x+r\left( 1-\frac{1}{c} \right) \right)$.
Writing $\frac{x}{c}$ in place of $x$ yields
$f(x) = f \left( \frac{1}{c}(x-r) + r \right)$
for every $c \in (0,\infty)$.

Let $x_{0}<r$ be fixed and let $c$ range over
$(0,\infty)$, then $\frac{1}{c}(x_{0}-r)+r$ ranges over $(-\infty,r)$, so
$f(x)$ is constant for $x<r$. Similarly, $f(x)$ is also constant for $x>r$.
\ep

\section{Open questions}
\lab{s:open}

The most important open question is the following. By Theorem \ref{t:meas} the
difficulty is to handle rotations.

\bq
\lab{q:meas}
Is every vertically rigid Lebesgue (Baire)
measurable function of the form $a+bx$ or $a+be^{kx}$ $(a,b,k \in
\RR)$ almost everywhere (on a comeagre set)? Or is this conclusion
true at least for Borel measurable functions, or Baire class~1
functions, or functions with at least one point of continuity?
\eq

\br
It would be more natural to replace vertical rigidity by
\emph{almost} vertical rigidity. However, it is not clear how this
should be defined, as a set can have a measure zero projection on
one line and positive measure projection on another.
\er

\bq Let $f$ be a vertically rigid function and $c>0$ such that there
exists an isometry between $graph(f)$ and $graph(cf)$ that is not a
translation (or also not a reflection). Is then $f$ of the form
$a+bx$? Or is this true for Borel, Lebesgue, or Baire measurable
functions? And if we assume the same for \emph{every} isometry
between $graph(f)$ and $graph(cf)$? \eq

Perhaps the following question can be answered in the negative
by an easy transfinite recursion. A positive answer to the analytic (see
e.g.~\cite{Ke} for the definition of analytic sets) version would answer
Question \ref{q:meas} for Borel functions.

\bq
\lab{q:rig}
Let $I \su S^1$ be the open subarc of arclength $\frac\pi2$
connecting $(0,1)$ and $(1,0)$. For a rigid set $H$ can $H \cap I$
be anything else but $\eset$, a point, $I$ minus a point, or $I$?
How about analytic, Borel, or Lebesgue (Baire) measurable rigid
sets?
\eq

\bq
What is the role of the uncountable set $C \su (0,\infty)$ in the
results of this paper? When is it sufficient to assume that it is infinite,
dense, sufficiently large finite, or contains a $c \neq 1$?
\eq

\br Let $c_0 \neq 1$. It is easy to see that there exists a
continuous $f$ satisfying $c_0f(x) = f(x+1)$ for every $x$. Indeed,
if we define $f$ to be an arbitrary continuous function on $[0,1]$
satisfying $c_0f(0) = f(1)$ then $f$ extends to $\RR$ in a unique
manner. Then $f$ is vertically rigid via horizontal translations for
the set $C = \{c_0^n : n\in\ZZ\}$. Hence it is not sufficient to
assume for Jankovi\'c's conjecture that $C$ is infinite.

There also exists a continuous nonlinear function $f$ whose graph
consists of two half lines starting from the origin so that
$graph(2f)$ is a rotated copy of $graph(f)$. \er

\bq
Is every horizontally rigid function of the form $a+bx$ or of the form
described in Theorem \ref{t:horiz}? Or is this true if we assume Borel,
Lebesgue, or Baire measurability? Is every continuous horizontally rigid
function of the form $a+bx$?
\eq

\bq
What can we say in higher dimensions?
\eq


\noindent \textbf{Acknowledgement} In the first place, we would like
to thank the anonymous referee for his valuable comments. His
suggestions vastly improved the paper. We are also indebted to
M.~Laczkovich for some helpful remarks.


\end{document}